\documentclass[12pt]{article}
\usepackage{amsfonts,amsmath,amsxtra,euscript}
\usepackage[english]{babel}
\author{Gleb G. Gusev\thanks{Partially supported by the grants
RFBR 10-01-00678, RFBR 08-01-00110-a, RFBR and SU HSE 09-01-12185-ofi-m, NSh-8462.2010.1.}}
\title{Euler characteristic of the bifurcation set for a polynomial of degree 2 or 3}
\date{}

\newtheorem{thm}{Theorem}

\newtheorem{prop}{Proposition}
\newtheorem{lem}{Lemma}

\newenvironment{dfntn} {\smallskip\noindent{\bf
Definition\/}.}{\smallskip\par}
\newenvironment{clm} {\smallskip\noindent{\bf
Claim\/}.}{\smallskip\par}
\newenvironment{exmpl} {\smallskip\noindent{\bf Example\/}.}{\smallskip\par}

\newenvironment{prf} {\noindent{\em Proof\/}.}{{ $\Box$}\smallskip\par}

\renewcommand{\a}{\alpha }

\renewcommand{\d}{\delta }
\newcommand{\D}{\Delta }

\renewcommand{\L}{\Lambda }

\newcommand{\s}{\sigma }

\newcommand{\p}{\pi}

\newcommand{\xx}{\mathbf x}
\newcommand{\zz}{\mathbf z}

\newcommand{\yy}{\mathbf y}

\newcommand{\kk}{\mathbf k}
\newcommand{\R}{\mathbb{R}}

\newcommand{\C}{\mathbb{C}}
\newcommand{\N}{\mathbb{N}}

\newcommand{\DDD}{\mathfrak{D}}

\newcommand{\Z}{\mathbb{Z}}
\newcommand{\ZZ}{\EuScript{Z}}
\newcommand{\WW}{\widetilde{W}}

\newcommand{\Vol}{\mathop{\mathrm{Vol}}\nolimits}

\begin{document}
\maketitle

\abstract{Assume that the coefficients of a polynomial in a complex variable are Laurent polynomials in some complex parameters. The parameter space (a complex torus) splits into strata corresponding to different combinations of coincidence of the roots of the polynomial. For generic Laurent polynomials with fixed Newton polyhedra the Euler characteristics of these strata are also fixed. We provide explicit formulae for the Euler characteristics of the strata in terms of the polyhedra of the Laurent polynomials in the cases of degrees 2 and 3. We also obtain some corollaries in combinatorial geometry, which follows from two different ways of computing the Euler characteristic of the bifurcation set for a reduced polynomial of degree 2.}

\section{Introduction}

Assume that the coefficients of a polynomial $\,P_{\zz}(t) = p_0(\zz) t^k + p_1(\zz) t^{k-1} + \ldots + p_k(\zz)\,$ are Laurent polynomials in the complex variables $\,(z_1, z_2, \ldots, z_n) = \zz\,$. The parameter space $(\C^*)^n$ ($\,\C^* = \C \setminus
0\,$) splits into strata corresponding to different combinations of coincidence of the roots of the polynomial $\,P_{\zz}$. For generic polynomials $\,p_0, p_1,
\ldots, p_k\,$ with fixed Newton polyhedra the Euler characteristics of these strata are also fixed. We provide explicit formulae for the Euler characteristics of the strata in terms of the polyhedra of the polynomials $p_i$ in the cases $k=2,
3$. This paper is an extension of \cite{euler} where the case $k=2$ was studied. We remind some usual notions first.

\begin{dfntn}
For a Laurent polynomial $\,S = \sum_{\kk\in \Z^n} s_\kk \zz^\kk,\,$ its Newton polyhedron $\D_S$ is the convex hull $\,\left <\, \{\kk \in \Z^n \mid s_\kk \neq 0\}\, \right > \subset \R^n\,$ of the set of integer points corresponding to nonzero coefficients.
\end{dfntn}

Denote by $\,\d_0, \d_1, \ldots, \d_k\,$ the Newton polyhedra of the polynomials $\,p_1, p_2,\ldots, p_k\,$ respectively. Consider the (Laurent) polynomial $P$ in $(n+1)$ variables that is defined by $\,P(\zz, t) = P_{\zz}(t).\,$ Denote by $\,\D \subset \R^{n+1}\,$ the Newton polyhedron of the polynomial $P$.

Denote by $\ZZ^n$ the set of primitive integer covectors in the dual space ${(\R^n)}^*$. For a polyhedron $\D$ and a covector $\,\a \in \ZZ^n,\,$ define the face $\D^{\a}$ of the polyhedron as the subset of $\D$ where $\a|_{\D}$ reaches its minimal value: $\,\D^{\a} = \{\xx\in \D \mid \a(\xx)=\min(\a|_{\D})\}.\,$ For a Laurent polynomial $\,S = \sum_{\kk\in \Z^n} s_\kk \zz^\kk,\,$ denote by $S^\a$ the (Laurent) polynomial $\,\sum_{\kk\in \D^\a_S} s_\kk \zz^\kk$.

\begin{dfntn}\label{non-deg}
A system of Laurent polynomials $\,P_1,P_2,\ldots,
P_k\,$ in $n$ complex variables is non-degenerate with respect to its Newton polyhedra $\,\D_{P_1}, \D_{P_2},\ldots, \D_{P_k}\,$ if for each covector $\,\a\in \ZZ^n\,$ and for each point $\,\zz \in (\C^*)^n\,$ such that $\,P_1^\a(\zz)=P_2^\a(\zz)=\ldots=P_k^\a(\zz)=0\,$ the system of covectors $\,dP_i^\a(\zz),\,$ $i=1,2,\ldots, k,\,$ is linear-independent.
\end{dfntn}

Let $\,S_1, S_2,\ldots, S_n\subset \R^n\,$ be a set of convex bodies. The Minkovskian sum of two bodies $S_1, S_2$ is defined by $\,S_1+S_2 := \{\xx_1 + \xx_2 \mid \xx_i \in S_i,\,i=1,2\}.\,$ The Minkovskian mixed volume of the bodies $\,S_1, S_2,\ldots, S_n\,$ is
\begin{multline*}
S_1 S_2\cdots S_n = \frac{1}{n!} [ \Vol_n(S_1+\ldots +S_n) - \sum_{i_1<\ldots < i_{n-1}} \Vol_n(S_{i_1}+\ldots +S_{i_{n-1}}) +\\
+ \sum_{i_1<\ldots < i_{n-2}} \Vol_n(S_{i_1}+\ldots +S_{i_{n-2}}) - \ldots + (-1)^{n-1} \sum_{i=1}^n \Vol_n(S_i)],
\end{multline*}
where $\Vol_n$ stands for the usual volume in $\R^n$. For a homogenous polynomial $\,T(x_1, x_2, \ldots, x_k) = \sum
\a_{i_1 i_2\ldots i_n}\,x_{i_1} x_{i_2} \cdots x_{i_n},\,$ of degree
$n$ we define $\,T(S_1, S_2, \ldots, S_k)$ as
$\sum \a_{i_1 i_2\ldots i_n}\,S_{i_1} S_{i_2} \cdots S_{i_n}$.

Consider a system of Laurent polynomials $\,P_1,P_2,\ldots, P_k\,$ in $n$ variables and the set $V = \{\zz\in (\C^*)^n \mid P_i(\zz)=0, i=1,2,\ldots, n\}$. The celebrating theorem of A.~Khovanskii (\cite{Khov}) claims the following. If the system $P_1,P_2,\ldots, P_k$ is non-degenerate with respect to its Newton polyhedra $\,\D_1,\D_2,\ldots, \D_k\,$ then the Euler characteristic of $V$ is
\begin{equation}
\label{chi} \chi(V) = n! (-1)^{n-k} Q^n_k(\D_1,\D_2,\ldots, \D_k),
\end{equation}
where $Q^n_k(x_1,x_2,\ldots, x_k)$\label{QRS} is the  the homogenous part of degree $n$ of the series $\prod_{i=1}^k \frac{x_i}{1-
x_i}.\,$ In particular,
\begin{equation}
\begin{split}
&Q^n_2(x,y) = \sum_{i=1}^{n-1} x^i
y^{n-i},\\
&Q^n_3(x,y,z) = \sum_{i,j,k \geq 1,\,i+j+k=n} x^i y^j
z^k,\\
&Q^n_4(x,y,z,t) = \sum_{i,j,k,l \geq 1,\,i+j+k+l=n} x^i
y^j z^k t^l.\\
\end{split}\notag
\end{equation}

\newpage
\section{The case of a polynomial of degree two}

In the case the degree of the polynomial $P$ is $\,k=2\,$ the parameter space $(\C^*)^n$ splits into $5$
strata --- $K$: $\deg(P_{\zz}) = 2$ and the roots of the polynomial $P_{\zz}$ are distinct; $L$: $\deg(P_{\zz}) = 2$ and the roots coincide; $M$: $\deg(P_{\zz}) = 1$; $N$: $\deg(P_{\zz}) = 0$; $O$: $P_{\zz}
\equiv 0$.

\begin{thm} \label{1Thm1}
For generic Laurent polynomials $p_i$ with fixed Newton polyhedra $\d_i$, $i=0,1,2$, the following equations hold:
\begin{equation}
\begin{split}
&\chi(K) = (-1)^n n!\,[\d_0^n + 2\d_*^n + Q^n_2(\d_0, \d_*)+
Q^n_2(\d_*, \d_2)+ Q^n_2(\d_0, \d_1) + Q^n_3(\d_0, \d_1, \d_2)],\\
&\chi(L) = (-1)^{n-1}n!\, [2\d_*^n + Q^n_2(\d_0, \d_*)+
Q^n_2(\d_*, \d_2)+ Q^n_2(\d_0, \d_1) + Q^n_3(\d_0, \d_1, \d_2)],\\
&\chi(M) = (-1)^{n-1}n!\, [\d_0^n + Q^n_2(\d_0, \d_1)],\\
&\chi(N) = (-1)^n n!\, [Q^n_2(\d_0, \d_1)+Q^n_3(\d_0, \d_1, \d_2)],\\
&\chi(O) = (-1)^{n-1}n!\,Q^n_3(\d_0, \d_1, \d_2),
\end{split} \notag
\end{equation} where $\,\d_* = \langle \d_1 \cup 1/2(\d_0 + \d_2)
\rangle$, $\langle \cdot \rangle\,$ denotes the convex hull,
$+$ stands for the Minkovskian sum, $\,Q^n_k(x_1,x_2,\ldots, x_k) = \left [\prod_{i=1}^k \frac{x_i}{1-
x_i}\right ]_n,\,$ $[\cdot ]_n$ is the the homogenous part of degree $n$ of the series under consideration.

A sufficient generality condition of the system of Laurent polynomials $p_i$ consists of the following sentences.

\begin{enumerate}
\item The systems of Laurent polynomials
$$
\{p_0,p_1,p_2\},\,\{p_0,p_1\},\,\{p_0\}
$$
are non-degenerate with respect to the systems of Newton polyhedra
$$
\{\d_0, \d_1, \d_2\}, \{\d_0, \d_1\} \mbox{, and } \{\d_0\}
$$
respectively.

\item The Laurent polynomial $P$ is non-degenerate with respect to its Newton polyhedron~$\D$.
\end{enumerate}
\end{thm}

\begin{prf}
The sets $\,O,\,N\sqcup O,\,\,M\sqcup N\sqcup O\,$ correspond to the systems of polynomials $\,\{p_0,p_1,p_2\},\,\{p_0,p_1\},\,\{p_0\}\,$ respectively, which are non-degenerate by the assumption of the theorem. We apply the equation (\ref{chi}) to these systems, use the additivity of the Euler characteristic and obtain $\chi(M)$, $\chi(N)$, $\chi(O)$.
The following idea is an analogue of the proof of lemma 2 from
\cite{Esterov} and provides $\,\chi(K), \chi(L)$.

Let $X$ be the subset of $\,(\C^*)^n\times \C\,$ defined by the equation $P=0$. Consider the projection $\,p \colon (\C^*)^n\times \C \to (\C^*)^n\,$ onto the first factor. Denote by $\p = p|_X$ its restriction to $X$. The Euler characteristics of the strata $\,K,L,M,N,O\,$ and the one of $X$ are related by two following linear equations. The first one

\begin{equation}\label{1eq1}
\chi(K) + \chi(L) + \chi(M) + \chi(N) + \chi(O) = 0
\end{equation}
follows from the additivity of the Euler characteristic (the left-hand side of (\ref{1eq1}) is the Euler characteristic of the torus $(\C^*)^n$).
One obtains the second equation computing the integral with respect to the Euler characteristic:
$$
\chi(X) = \int_{(\C^*)^n} \chi(\p^{-1}(\zz))\,d\chi.
$$
The pre-image of a point $\,\zz \in K\,$ consists of two points, the one of a point $\,\zz \in L\sqcup M\,$ consists of one point, the  pre-image of $N$ is empty, and finally, for a point $\,\zz \in O,\,$ one has: $\,\p^{-1}(\zz) \cong \C\,$ and $\,\chi(\p^{-1}(\zz)) = 1.\,$ Therefore, one obtains:
\begin{equation}\label{1eq2}
\chi(X) = 2\chi(K)+ \chi(L)+ \chi(M)+ \chi(O).
\end{equation}
Combining the linear equations (\ref{1eq1}), (\ref{1eq2}) one obtains:
\begin{equation}
\label{1chi}\chi(L) = -\chi(X) - \chi(M) - 2\chi(N) - \chi(O).
\end{equation}

Let us find $\chi(X)$. The set $X$ consists of two strata $\,X_1 =
X \setminus \{(\zz, t)\mid t=0\}\,$ and $\,X_2 = X \cap \{(\zz,
t)\mid t=0\}\,$. One has:
\begin{equation}\label{1X}
\chi(X) = \chi(X_1) + \chi(X_2).
\end{equation}
The set $X_2$ is given in $(\C^*)^n$ by the non-degenerate equation $\,p_2 = 0,\,$ therefore applying the equation~(\ref{chi}) one obtains
\begin{equation}\label{1chi_X2}
\chi(X_2) = (-1)^{n-1}n!\,\d_2^n.
\end{equation}
The stratum $X_1$ is given in the torus $(\C^*)^{(n+1)}$ by the equation $\,P=0.\,$ Applying the equation~(\ref{chi}) to the polynomial $P$, one obtains
\begin{equation}
\label{1chi_X1}\chi (X_1) = (-1)^{n}(n+1)!\,\D^{n+1}.
\end{equation}

Let us express the volume of $\D \subset \R^{n+1}$ in terms of the polyhedra $\d_i$. Denote by~$k_t$ the coordinate in the space $\R^{n+1}$ that corresponds to the variable $t$. Denote by~$v_t$ the vector in $\R^{n+1}$ that has only one non-zero coordinate $k_t = 1$. The polyhedra $\d_i$ lie in the hyperplane $k_t=0$, and the polyhedron $\D$ is the convex hull of their parallel shifts along the $k_t$-axis: $\,\D = \langle\, (\d_0+ 2v_t) \cup (\d_1+ v_t) \cup \d_2 \,\rangle\,$.

\begin{clm}
The intersection of the hyperplane $\{k_t=1\}$ and the polyhedron $\D$ is the polyhedron $\,\d_*+ v_t$, where $\d_* = \langle \d_1 \cup
1/2(\d_0 + \d_2) \rangle$.
\end{clm}

\begin{prf}
An arbitrary point $\,(\kk,\a)\in \D\,$ can be expressed as $(\kk,\a) = \a_0(\kk_0,2)+ \a_1(\kk_1,1)+\a_2(\kk_2,0),\,$ where $\kk_i\in \d_i$, $\a_i>0,\,$
$\,i=1,2,3,\,$ and $\sum \a_i = 1$. Assume that the point
$(\kk,\a)$ lie in the hyperplane $\{k_t=1\}$, other words, $\a=1$. Then one has $\,\a_0=\a_2\,$ and therefore the point $\,(\kk,\a)\,$ is a convex combination of the points $\,\frac{1}{2}((\kk_0,2)+ (\kk_2,0))\,$ and
$(\kk_1,1)$ of the polyhedra $\frac{1}{2}(\d_0 + \d_2)+ v_t$ and $\d_1 + v_t$ respectively.
\end{prf}

The volumes of the parts $\,\langle\, (\d_0+ 2v_t) \cup (\d_*+ v_t)
\,\rangle\,$ and $\,\langle\, (\d_*+ v_t) \cup \d_2 \,\rangle,\,$ which form the polyhedron $\D$, can be obtained by the following formula.

\begin{lem}\label{1lem}
Suppose that polyhedra $\L_0,\,\L_1\subset \R^{n+1}$ lie in the $n$-dimensional hyperplanes given in the space $\R^{n+1}$ by the equations $x_1
=0,\,x_1=1$ respectively. Then the $(n+1)$-dimensional volume of their convex hull $\L$ equals
\begin{equation}
\label{1V} \L^{n+1}=\frac{1}{n+1}\,(\L_0^n+ Q^n_2(\L_0, \L_1)+
\L_1^n).
\end{equation}
\end{lem}

\begin{prf}
For $\a \in [0,1],$ the section $\L_{\a}$ of the polyhedron $\L$
by the hyperplane $\,x_1=\a\,$ consists of points of the form $\,(1-\a)\xx+ \a\yy,\,$ where $\,\xx \in \L_0,\,$ and $\,\yy \in \L_1\,$. It follows that
$\,\L_{\a}=(1-\a)\L_0 + \a \L_1,\,$ and therefore
$$
\L^{n+1}=\,\int_{[0,1]} \,(\L_{\alpha})^n d\alpha
=\,\int_{[0,1]}\left [\sum_{i=0}^n C_n^i \a ^i (1-\a
)^{n-i}(\L_1)^i(\L_0)^{n-i} d\a \right ].
$$
Taking into account that
$$
\int_{[0,1]}(\a^i (1-\a)^{n-i})d\a=
i!\,(n-i)!/(n+1)!=1/((n+1)C_n^i),
$$
one obtains the formula (\ref{1V}).
\end{prf}

Applying (\ref{1V}) to the parts of the polyhedron $\D$ and using the invariance of the mixed volume with respect to parallel translations of its arguments one obtains
\begin{equation}\label{1del}
\D^{n+1} = 1/(n+1)\,(\d_0^n + 2\d_*^n + \d_2^n + Q^n_2(\d_0, \d_*)
+ Q^n_2(\d_*, \d_2)).
\end{equation}
The equations of Theorem~\ref{1Thm1} for $\,\chi(K), \chi(L)\,$ follows now from (\ref{1eq1}), (\ref{1chi}), (\ref{1X}), (\ref{1chi_X2}),
(\ref{1chi_X1}), (\ref{1del}).
\end{prf}

\subsection{The corollaries in combinatorial geometry.}

For a reduced polynomial of degree two $\,P_{\zz}(t) = t^2 +
p_1(\zz)\,t + p_2(\zz)\,$ the strata $\,M,N,O\subset (\C^*)^n\,$ are empty and thus
their Euler characteristics are equal to zero. This fact respects the formulae of Theorem \ref{1Thm1} (taking into account that
$\,\d_0=\{0\}$). The stratum $L$ is given by the equation ${\,p_1^2 - 4 p_2
= 0}.\,$ It follows from Theorem~\ref{1Thm1} that for generic $\,p_1, p_2\,$ the Euler characteristic of $L$ is given by the formula:
\begin{equation}\label{1chi_L_1}
\chi(L) = (-1)^{n-1}n!\, [2\d_*^n + Q^n_2(\d_*, \d_2)],
\end{equation}
where $\,\d_* = \langle \d_1 \cup 1/2\d_2 \rangle$.

Using the method of toroidal compactifications that was provided by A. Khovanskii for studying the invariants of non-degenerate hypersurfaces in the torus $(\C^*)^n$ (see~\cite{Khov}), one can prove another formula for the Euler characteristic of the stratum $L$ by induction on $n$:
\begin{equation}\label{1chi_L_2}
\chi(L) = (-1)^{n-1}n!\, [(2\d_*)^n - Q^n_2(2\d_*, \d_1) +
Q^n_2(\d_1, \d_2)].
\end{equation}

The coexistence of the two formulae (\ref{1chi_L_1}), (\ref{1chi_L_2}) does not induce contradiction because the polyhedra $\,\d_*, \d_1,
\d_2\,$ are related.

\begin{prop}
Consider a pair of arbitrary convex bodies $\,S_1,S_2\subset \R^n\,$ and the convex hull of their union $\,S_0 = \langle S_1 \cup S_2 \rangle.\,$ Then one has
\begin{equation}\label{1kombi}
R^n(S_0, S_1, S_2) = 0,
\end{equation}
where $\,R^n(x_0, x_1, x_2) = (2^n -2) x_0^n + Q^n_2(x_1, 2x_2) -
Q^n_2(2x_0, x_1) - Q^n_2(x_0, 2x_2)$.
\end{prop}

\begin{prf}
Assume first, that $\,S_1, S_2\,$ are integer polyhedra. Consider a generic pair of Laurent polynomials $\,p_1, p_2\,$ with fixed Newton polyhedra
$\,S_1, 2S_2\,$ respectively. There are the formulae (\ref{1chi_L_1}), (\ref{1chi_L_2}) for the Euler characteristic of the stratum $\,{L\subset (\C^*)^n},\,$ that corresponds to the coincidence of the roots of the polynomial
$\,P_{\zz}(t) = t^2 + p_1(\zz)\,t + p_2(\zz),\,$ and these formulae imply the equation (\ref{1kombi}).

Assume now, that $\,S_1, S_2\,$ are polyhedra with rational coordinates of all their vertices. There exists a natural number $k \in \N$ such that $\,kS_1, kS_2\,$ are integer polyhedra. One has: $\,R^n(S_0, S_1, S_2) =
(1/k^n) R^n(kS_0, kS_1, kS_2) = 0.\,$

Finally, let's consider the general case of an arbitrary pair $\,S_1, S_2\,$ of convex bodies. Consider a couple of consequences   $\,(S_1^i),\,(S_2^i)\,$ of polyhedra with rational coordinates of all their vertices such that $\,S_1^i
\xrightarrow[i \rightarrow \infty]{} S_1,\, S_2^i \xrightarrow[i
\rightarrow \infty]{} S_2\,$ (one can choose the n-dimensional volume of symmetric difference of two bodies as the metric on the set of convex bodies in the space $\R^n$). Let $\,S_0^i = \langle S_1^i \cup S_2^i \rangle\,$ Then one has: $\,R^n(S_0, S_1, S_2) =
\lim_{i\to \infty} R^n(S_0^i, S_1^i, S_2^i) = 0.\,$
\end{prf}

\begin{exmpl}
Let  $n=2$. Then $\,R^n(x_0, x_1, x_2) = 2x_0^2 + 2x_1x_2 -
2x_0x_1 - 2x_0x_2 = 2(x_0 - x_1)(x_0 - x_2).\,$ Thus one obtains the following corollary: for any three convex figures $\,S_0,S_1,S_2\subset \R^2\,$ that are connected by the relation $\,S_0 = \langle S_1 \cup S_2
\rangle\,$ one has:
$$
(S_0-S_1)(S_0-S_2) = 0.
$$
\end{exmpl}

\section{The case of a polynomial of degree three}

For $\,k=3,\,$ the parameter space $(\C^*)^n$ splits into $8$
strata. The strata $\,K,L,M,N,O\,$ are defined in the previous section, and one has 3 new strata in addition --- $H$: $\deg(P_{\zz}) = 3$ and the roots of the polynomial $P_{\zz}$ are distinct; $I$: $\deg(P_{\zz}) = 3$ and the polynomial has a double root; $J$: $\deg(P_{\zz}) = 3$ and the polynomial has a triple root.

Consider the inclusion $\,\R^n\subset \R^{n+1}\subset \R^{n+2},\,$ where the first space is equipped by the coordinates $\,\kk =
(k_1,k_2\ldots,k_n),\,$ the second one is equipped by the additional coordinate $k_t$ (see the previous section), and the third one has moreover the additional coordinate $k_\s$ that corresponds to an additional complex variable $\s$.

Denote by~$v_t$ the vector in $\R^{n+2}$ that has only one non-zero coordinate $k_t = 1$. Denote by $v_\s$ the vector in $\R^{n+2}$ that has only one non-zero coordinate~${v_\s=1}$. Denote by $\D_{1,2,3}$ the Newton polyhedron of the Laurent polynomial $\,{p_1 t^2
+p_2 t+ p_3}.\,$ Denote by $\DDD_i$ the Newton polyhedron of the Laurent polynomial $\,\s+p_i
t^{3-i}\,$ ($\,i=0,1,2,3\,$) in $n+2$ variables.  One has: $\D_{1,2,3} = \left<(\d_1+2v_t)\cup (\d_2+ v_t) \cup \d_3\right>$, $\,\DDD_i = \left< v_\s \cup (\d_i + (3-i)v_t) \right>$.

\begin{thm}
For generic Laurent polynomials $p_i$ with fixed Newton polyhedra $\d_i$, $i=0,1,2,3$, one has:
\begin{equation*}
\begin{split}
\chi(H) = &(-1)^n n!\,[(n+1)(n+2) Q^{n+2}_4(\DDD_0, \DDD_1,
\DDD_2,
\DDD_3)+\\
          &+(n+1)(\D^{n+1}+ Q^{n+1}_2(\d_0, \D_{123}))-\\
          &-2\d_0^n - \d_3^n - Q^n_2(\d_0, \d_3) -
Q^n_3(\d_1,\d_2,\d_3) - Q^n_4(\d_0,\d_1, \d_2, \d_3)],\\
\chi(I) = &(-1)^{n-1} n!\, [2(n+1)(n+2) Q^{n+2}_4(\DDD_0, \DDD_1,
\DDD_2, \DDD_3)+\\
          &+(n+1)(\D^{n+1} + Q^{n+1}_2(\d_0, \D_{123}))-\\
          &-3\d_0^n -\d_3^n - Q^n_2(\d_0,\d_3)- 2Q^n_3(\d_1,\d_2,\d_3)-
2Q^n_4(\d_0,\d_1,\d_2,\d_3)],\\
\chi(J) = &(-1)^n n!\,[(n+1)(n+2) Q^{n+2}_4(\DDD_0, \DDD_1,
\DDD_2,
\DDD_3)- Q^n_3(\d_1,\d_2,\d_3)-\\
          &- Q^n_4(\d_0,\d_1, \d_2, \d_3)],\\
\chi(K) = &(-1)^{n-1} n!\, [(n+1)Q^{n+1}_2(\d_0, \D_{123}) -
\d_0^n
- Q^n_2(\d_0, \d_3)+\\
          &+Q^n_3(\d_0,\d_1,\d_2) + Q^n_4(\d_0,\d_1, \d_2,
\d_3)],\\
\chi(L) = &(-1)^n n!\, [(n+1)Q^{n+1}_2(\d_0, \D_{123}) -
2\d_0^n - Q^n_2(\d_0, \d_3)-\\
          &- Q^n_2(\d_0, \d_1) + Q^n_3(\d_0,\d_1,\d_2) +
Q^n_4(\d_0,\d_1, \d_2, \d_3)],\\
\chi(M) = &(-1)^n n!\, [Q^n_2(\d_0, \d_1) + Q^n_3(\d_0,\d_1,\d_2)],\\
\chi(N) = &(-1)^{n-1} n!\,[Q^n_3(\d_0,\d_1,\d_2) +
Q^n_4(\d_0,\d_1, \d_2, \d_3) ],\\
\chi(O) = &(-1)^n n!\,Q^n_4(\d_0,\d_1, \d_2, \d_3),\\
\end{split}
\end{equation*} where the polynomials $\,Q^n_i\,$ are defined on the page~\pageref{QRS}.

A sufficient generality condition consists of the following requirements. The systems of Laurent polynomials in $n$ variables
$$
\{p_0,p_1,p_2,p_3\},\,\{p_0,p_1,p_2\},\,\{p_1,p_2,
p_3\}\,\{p_0,p_1\}\,\{p_0,p_3\}\,\{p_0\},\{p_3\}
$$
are non-degenerate with respect to its systems of Newton polyhedra
$$
\{\d_0, \d_1, \d_2, \d_3\}, \,\{\d_0, \d_1, \d_2\}, \,\{\d_1,
\d_2, \d_3\}, \,\{\d_0, \d_1\}, \,\{\d_0, \d_3\}, \,\{\d_0\},
\,\{\d_3\}
$$
respectively. The systems of Laurent polynomials in $n+1$ variables
$$
\{p_0, p_1 t^2 +p_2 t+ p_3\}, \,\{P\}
$$
are non-degenerate with respect to its systems of Newton polyhedra
$$
\{\d_0, \D_{123}\},\, \{\D\}
$$
respectively. Finally, the system of Laurent polynomials in ${n+2}$ variables $\,\{\s- 3p_0 t^3, \s+ p_1 t^2, \s- p_2 t, \s+ 3p_3
\}\,$ is non-degenerate with respect to its Newton polyhedra $\,\{\DDD_0, \DDD_1, \DDD_2, \DDD_3\}.$
\end{thm}

\begin{prf}
The equations of the Theorem~\ref{1Thm1} appear as the solution of a system of 8 independent linear equations. One obtains the first of them calculating the Euler characteristic of the torus $(\C^*)^n$:
\begin{equation}\label{1equat1}
\chi(H)+ \chi(I)+ \chi(J)+ \chi(K) + \chi(L) + \chi(M) + \chi(N) +
\chi(O) = 0
\end{equation}

The sets $\,O,\,N\sqcup O,\,\,M\sqcup N\sqcup O,\,\,K\sqcup
L\sqcup M\sqcup N\sqcup O\,$ are given by the systems of Newton polyhedra $\,\{p_0,p_1,p_2,p_3\},\,\{p_0,p_1,p_2\},$
$\{p_0,p_1\},\,$ and $\,\{p_0\}\,$ respectively. Applying the equation (\ref{chi}) to these systems one obtains 4 equations:
\begin{equation}\label{1equat2}
\begin{split}
&\chi(O) = (-1)^n n!\,Q^n_4(\d_0, \d_1, \d_2, \d_3),\\
&\chi(N) + \chi(O) = (-1)^{n-1} n!\,Q^n_3(\d_0, \d_1, \d_2),\\
&\chi(M) + \chi(N) + \chi(O) = (-1)^n n!\,Q^n_2(\d_0, \d_1),\\
&\chi(K) + \chi(L) + \chi(M) + \chi(N) + \chi(O) = (-1)^{n-1} n!\,
\d_0^n.
\end{split}
\end{equation}

One obtains two extra equations calculating the following integrals with respect to the Euler characteristic. Consider the subsets $\,Y, Z \subset (\C^*)^n\times \C,\,$ where $Y =
\{(\zz, t)\mid \deg P_\zz \leq 2,\, P_\zz(t) =0\},\,$ $Z = \{(\zz,
t)\mid P_\zz(t) =0\}$. Let $\,\p_1\colon Y \to
(\C^*)^n,\,\p_2\colon Z \to (\C^*)^n\,$ be the restrictions of the projection
$\p\colon (\C^*)^n\times \C \to (\C^*)^n$ onto the first factor.
Let us calculate the integrals with respect to the Euler characteristic:
$$
\chi(Y) = \int_{(\C^*)^n} \chi(\p_1^{-1}(\zz))\,d\chi,\quad
\chi(Z) = \int_{(\C^*)^n} \chi(\p_2^{-1}(\zz))\,d\chi.
$$
For $\,\zz \in H\sqcup I\sqcup J,\,$ one has: $\,\p_1^{-1}(\zz) =
\emptyset.\,$ The pre-image $\,\p_2^{-1}(\zz)\,$ consists of three points for $\,\zz \in H,\,$ consists of two points for $\,\zz \in I,\,$ and consists of one point for
$\,\zz \in J.\,$ For $\,\zz \in K,\,$ the pre-images under both maps consist of two points; for $\,\zz \in L\sqcup
M,\,$ consist of one point; for $\,\zz \in N,\,$ the pre-images are empty; finally, for $\,\zz \in O,\,$ the pre-images are isomorphic to $\C$, thus
$\chi(\p_i^{-1}(\zz)) = 1,\,i=1,2$. Therefore, one has:
\begin{equation}\label{1integ}
\begin{split}
&\chi(Y) = \chi(O) + \chi(M) + \chi(L)+ 2\chi(K),\\
&\chi(Z) = \chi(O) + \chi(M) + \chi(L)+ 2\chi(K)+ \chi(J)+
2\chi(I) + 3\chi(H).
\end{split}
\end{equation}

Consider the decompositions $\,Y=Y_1\sqcup Y_2,\,Z=Z_1\sqcup Z_2,\,$ where $\,{Y_1 = Y\cap (\C^*)^{n+1}}$, $\,Z_1 = Z\cap (\C^*)^{n+1},\,Y_2 =
\{(\zz,t)\in Y\mid t=0\},\,Z_2 = \{(\zz,t)\in Z\mid t=0\}.\,$
The strata $\,Y_1, Z_1\,$ are given in $\,(\C^*)^{n+1}\,$ by the systems $\,\{p_0, p_1 t^2 + p_2 t + p_3\},\,\{P\}\,$
respectively. The strata $\,Y_2, Z_2\,$ are given in $\,(\C^*)^n\,$ by the systems $\,\{p_0,p_3\},\,$ and $\,\{p_3\}\,$ respectively. Applying the equations (\ref{chi}) one obtains the Euler characteristics of the strata $\,Y,Z.\,$ One substitutes the answers into the formulae (\ref{1integ}) and obtains the two linear equations:
\begin{multline}\label{1equat3}
\chi(O) + \chi(M) + \chi(L)+ 2\chi(K) =\\
= (-1)^{n-1} n!\,((n+1) Q^{n+1}_2(\d_0, \D_{123})- Q^n_2(\d_0,
\d_3)),
\end{multline}
\begin{multline}\label{1equat4}
\chi(O) + \chi(M) + \chi(L)+ 2\chi(K)+ \chi(J)+ 2\chi(I) +
3\chi(H) =\\
= (-1)^n n!\,((n+1)\D^{n+1} - \d_3^n).
\end{multline}
The equations (\ref{1equat2}), (\ref{1equat3}) provide the characteristics of the strata $\,K, L, M, N, O.\,$ One obtains the characteristics of the strata $\,H,I,J\,$ using the equations (\ref{1equat1}), (\ref{1equat4}), and the following final (eigth)  equation.

Consider the set
$$
W = \{(\zz,t_0)\in (\C^*)^n\times \C \mid P_\zz (t) = (\partial P_\zz/\partial t)(t_0) = (\partial^2 P_\zz/\partial t^2)(t_0)= 0\}.
$$
Arguments similar to the calculation of $\chi(Y)$ and $\chi(Z)$ provide that
\begin{equation}\label{1chi_W}
\chi(W) = \chi(O) + \chi(J).
\end{equation}
One has: $\,W=W_1\sqcup W_2,\,$ where $\,W_1 = W\cap
(\C^*)^{n+1},\,W_2 = \{(\zz,t)\in W\mid t=0\}.\,$ The stratum $W_2$
is given in $\,(\C^*)^n\,$ by the system $\,\{p_1,p_2,p_3\}.\,$ Therefore,
\begin{equation}\label{1W_2}
\chi(W_2) = (-1)^{n-1}n!\,Q^n_3(\d_1,\d_2,\d_3).
\end{equation}

Let us find $\chi(W_1)$. The set $W_1$ is given in the torus
$(\C^*)^{n+1}$ by the system $\,\{P = \partial P/\partial t =
\partial^2 P/\partial t^2 = 0\}.\,$ This system of Laurent polynomials may be degenerate with respect to its Newton polyhedra, thus one can not apply to it the equation~(\ref{chi}). One can see that the above system is equivalent to the following one, which may be still degenerate: $\,\{3p_0 t^3 + 3p_3 = 3p_0 t^3 - p_2 t = 3p_0 t^3+ p_1 t^2 = 0\}.\,$
Consider now the system of Laurent polynomials in $n+2$ variables
\begin{equation}\label{1system}
\{\s - 3p_0 t^3, \s + p_1 t^2 , \s - p_2 t, \s + 3p_3 \}.
\end{equation}

Denote the set of common zeroes of this system in the torus $(\C^*)^{n+2}$ by $\WW$. The projection $\,s \colon (\C^*)^{n+1}\times \C_\s \to (\C^*)^{n+1}\,$ onto the hyperplane of the coordinates $\,(\zz, t)\,$ embeds the set $\WW$ into the stratum $W_1$. Namely, let us consider the decomposition $\,W_1 = W_1'\sqcup W_1'',\,$
where $\,W_1' = \{(\zz, t)\in W_1\mid p_0(\zz)\neq 0\}\,$ and $\,W_1''= \{(\zz, t)\in W_1\mid p_0(\zz) = 0\}.\,$ Then one has: $\,s(\WW) = W_1'.\,$ The set $W_1''$ is given in $\C^{n+1}$ by the following system of equations in $n$ variables: $\,\{p_0=p_1=p_2=p_3=0\}.\,$ Therefore, $W_1''$ is a fibration with the fibre $\C^*$ and $\,\chi(W_1'')=
0.\,$ It follows that:
\begin{equation}\label{1WW}
\chi(W_1) = \chi(W_1') = \chi(\WW).
\end{equation}

\begin{clm}
For generic Laurent polynomials $p_i$ with fixed Newton polyhedra $\,\d_i,\,\,i=0,1,2,3,\,$ the system of Laurent polynomials (\ref{1system}) is non-degenerate with respect to its Newton polyhedra $\,\DDD_0,\DDD_1, \DDD_2, \DDD_3.$
\end{clm}

\begin{prf}
Denote by $L_i$ the affine space of Laurent polynomials $S$ in $n+1$ variables such that $\,\D_S \subset \DDD_i,\,i=0,1,2,3.\,$ The set of non-degenerate systems $\,\{S_0, S_1, S_2,S_3\}\,$ of Laurent polynomials with fixed Newton polyhedra $\,\DDD_0,\DDD_1,\DDD_2,\DDD_3\,$ form a Zariski open subset $U$ in the product $\Pi = \prod_{i=0}^3 P(L_i)$ of projectivizations
of the affine spaces. The formulae (\ref{1system}) provide an inclusion of the set $\Psi$ of systems $\,\{p_0, p_1, p_2, p_3\}\,$ of Laurent polynomials with fixed Newton polyhedra $\,\d_0, \d_1, \d_2, \d_3\,$ into the variety $\Pi$. This inclusion realizes $\Psi$ as an open subset of $\Pi$. The set of non-degenerate systems of the form (\ref{1system}) is $U\cap\Psi$, and therefore, is a Zariski open subset in $\Psi$.
\end{prf}

Applying the equation (\ref{chi}) to the system (\ref{1system}), and
taking into account (\ref{1chi_W}), (\ref{1W_2}), (\ref{1WW}), one gets the last needed linear equation:
\begin{multline}
\chi(O) + \chi(J) = (-1)^n n!\,[(n+1)(n+2)Q^{n+2}_4(\DDD_0,
\DDD_1,
\DDD_2, \DDD_3) - \\
 -Q^n_3(\d_1,\d_2,\d_3)].
\end{multline}
\end{prf}

\end{document}